\def\draft{n}
\documentclass[12pt]{amsart}
\usepackage{fullpage,amssymb,epic,eepic}

\theoremstyle{plain}

\newtheorem{proposition}{Proposition}[section]
\newtheorem{lemma}[proposition]{Lemma}

\theoremstyle{definition}

\newtheorem{question}{Question}

\theoremstyle{remark}

\newtheorem{remark}[proposition]{Remark}

\def\printname#1{
	\if\draft y
		\smash{\makebox[0pt]{\hspace{-0.5in}
			\raisebox{8pt}{\tt\tiny #1}}}
	\fi
}

\newlength{\standardunitlength}
\catcode`\@=11
\long\def\@makecaption#1#2{%
    \vskip 10pt
    \setbox\@tempboxa\hbox{
      \small\sf{\bfcaptionfont #1. }\ignorespaces #2}%
    \ifdim \wd\@tempboxa >\captionwidth {%
        \rightskip=\@captionmargin\leftskip=\@captionmargin
        \unhbox\@tempboxa\par}%
      \else
        \hbox to\hsize{\hfil\box\@tempboxa\hfil}%
    \fi}
\font\bfcaptionfont=cmssbx10 scaled \magstephalf
\newdimen\@captionmargin\@captionmargin=2\parindent
\newdimen\captionwidth\captionwidth=\hsize
\catcode`\@=12

\def\lbl#1{\label{#1}\printname{#1}}

\def\eqdef{\overset{\text{def}}{=}}

\def\eqdef{\overset{\text{def}}{=}}

\def\BZ{\mathbb Z}
\def\BQ{\mathbb Q}

\def\G{\mathcal G}

\def\T{\mathcal T}
\def\M{\mathcal M}
\def\F{\mathcal F}

\def\T{\mathcal T}

\def\Ga{\Gamma}
\def\S{\Sigma}
\def\La{\Lambda}

\def\b{\beta}
\def\e{\epsilon}
\def\a{\alpha}
\def\g{\gamma}

\def\ov#1{\overline{#1}}

\def\t{\tau}
\def\ihs{integral homology 3-sphere}
\def\fti{finite type invariants}
\def\la{\langle}
\def\ra{\rangle}
\def\LT{lantern identity}
\def\lt{lantern identity}

\def\FT#1{\F^T_{#1} \M}

\def\GT#1{\G^T_{#1} \M}

\begin{document}


\title[Applications of the lantern identity]{Applications of
 the lantern identity}

\author{Stavros Garoufalidis}
\address{Department of Mathematics \\
         Brandeis University \\
         Waltham, MA 02254-9110, U.S.A. }
\email{stavros@oscar.math.brandeis.edu}
\thanks{The  author was partially supported by NSF grant 
       DMS-95-05105.\newline
       This and related preprints can also be obtained at
{\tt http://www.math.brown.edu/$\sim$stavrosg } \newline
1991 {\em Mathematics Classification.} Primary 57N10. Secondary 57M25.
\newline
{\em Key words and phrases:} lantern identity, finite type invariants.
}


\date{This edition: \today ; First edition: August 1, 1997  }

\begin{abstract}
The purpose of this note is to unify the role of the lantern identity
in the proof of several finiteness theorems. In particular, we show 
that for every nonnegative integer $m$, the vector space
(over the rationals) of type $m$ (resp. $T$-type $m$) invariants of 
integral homology 3-spheres are finite dimensional. These results have already
been obtained by \cite{Oh} and  \cite{GL2} respectively; our derivation
however is simpler, conceptual and relates to
several other applications of the lantern identity.
\end{abstract}

\maketitle
\tableofcontents


\section{The lantern identity}
\lbl{sec.lantern} 

In his seminal work using the lantern identity, D. Johnson 
showed that the Torelli group is  finitely generated \cite{Jo4}
and further determined its abelianization \cite{Jo5}.
 
In the present  note we use the lantern identity in order to give a
simple proof of the following two {\em finiteness} 
theorems:
\begin{itemize}
\item
The vector space of type $m$ invariants of \ihs s is finite dimensional
(for every $m$), due to T. Ohtsuki \cite{Oh}.
\item
The vector space of $T$-type $m$ invariants of \ihs s is finite 
dimensional (for every $m$), due to the author and J. Levine, \cite{GL2}. 
\end{itemize}
Our proof relates D. Johnson's finiteness theorems to the above
mentioned  results and explains in a conceptual way the relation
of the lantern identity and finiteness theorems.
For completeness, in Section \ref{sec.comp} we  mention two more 
applications of the \lt\ to finiteness theorems.
We also offer exercises and hints for the reader to figure out  forms of a 
{\em fictional} \lt\ (which would still prove the various finiteness results),
hoping to bring together a rather diverse audience and unify the use of
the lantern identity in the proofs of finiteness theorems.

We begin by recalling the lantern identity, 
which is usually written as an identity in the group of 
{\em framed
pure $3$-strand braids} $P_3$ on the plane, or in the
group of  {\em framed pure $4$-strand braids} $P_4(S^2)$ on the $2$-sphere
(otherwise known as the mapping class group of a $2$-sphere
with $4$ boundary components).  With the notation
of  Figures \ref{lantern3}, \ref{lantern2} we have:
\begin{eqnarray*}
\t_{123}\t_1\t_2\t_3 & = & \t_{12} \t_{13} \t_{23} 
\hspace{1cm}\text{ in } P_3 \\
 \epsilon_4 \epsilon_1 \epsilon_2 \epsilon_3 & = & \alpha \beta \gamma
\hspace{1.7cm}\text{ in } P_4(S^2)    
\end{eqnarray*} 
Note that for an unoriented  simple closed curve $c$,
a right-handed (resp. left-handed) Dehn twist on $c$ represents an element of 
the mapping class group denoted by $c$ (resp. $c^{-1}$), and  
that we compose maps from left
to right (contrary to the usual way which composition of functions is
written) and braids from top to bottom.
There is a map
$P_{n+1} \to P_n(S^2)$ (which places the last strand at infinity)
which maps one version of the lantern identity on the other, as is obvious
from the left hand side of Figure \ref{lantern2}.

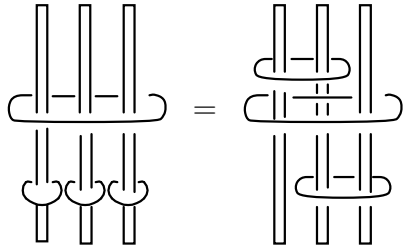
\begin{figure}[htpb]
$$ \printname{lantern3}
	\setlength{\unitlength}{0.03\standardunitlength}
	\begin{array}{c}  \hspace{-1.7mm}
        	\raisebox{-8pt}{\begingroup\makeatletter\ifx\SetFigFont\undefined
\def\x#1#2#3#4#5#6#7\relax{\def\x{#1#2#3#4#5#6}}%
\expandafter\x\fmtname xxxxxx\relax \def\y{splain}%
\ifx\x\y   
\gdef\SetFigFont#1#2#3{%
  \ifnum #1<17\tiny\else \ifnum #1<20\small\else
  \ifnum #1<24\normalsize\else \ifnum #1<29\large\else
  \ifnum #1<34\Large\else \ifnum #1<41\LARGE\else
     \huge\fi\fi\fi\fi\fi\fi
  \csname #3\endcsname}%
\else
\gdef\SetFigFont#1#2#3{\begingroup
  \count@#1\relax \ifnum 25<\count@\count@25\fi
  \def\x{\endgroup\@setsize\SetFigFont{#2pt}}%
  \expandafter\x
    \csname \romannumeral\the\count@ pt\expandafter\endcsname
    \csname @\romannumeral\the\count@ pt\endcsname
  \csname #3\endcsname}%
\fi
\fi\endgroup
\begin{picture}(5522,3369)(0,-10)
\thicklines
\path(3721,2412)(3721,3342)(3871,3342)(3871,2412)
\path(3721,2142)(3721,1752)
\path(3871,2112)(3871,1782)
\path(4321,2412)(4321,3342)(4471,3342)(4471,2412)
\path(4921,1842)(4921,3342)(5071,3342)(5071,1842)
\path(4321,1512)(4321,762)
\path(4471,1542)(4471,792)
\path(4921,1542)(4921,762)(4921,762)
\path(5071,1542)(5071,732)
\path(3721,1512)(3721,12)(3871,12)(3871,1542)
\path(4321,522)(4321,12)(4471,12)(4471,522)
\path(4921,522)(4921,12)(5071,12)(5071,522)
\path(3961,2592)(4261,2592)
\spline(3631,2082)
(3391,2082)(3301,1932)
	(3331,1752)(3421,1692)(5401,1692)
	(5491,1782)(5521,1902)(5491,2022)
	(5371,2082)(5281,2082)
\path(3991,2052)(4801,2052)
\path(4321,2232)(4321,2112)
\path(4471,2232)(4471,2112)
\path(4321,1962)(4321,1812)
\path(4471,1962)(4471,1812)
\spline(4261,942)
(4081,942)(3991,792)
	(4081,642)(5281,642)(5371,792)
	(5281,942)(5101,942)
\spline(3691,2592)
(3511,2592)(3421,2442)
	(3511,2292)(4711,2292)(4801,2442)
	(4711,2592)(4531,2592)
\path(4561,942)(4831,942)
\path(1621,1842)(1621,3342)(1771,3342)(1771,1842)
\path(1021,1512)(1021,762)
\path(1171,1542)(1171,792)
\path(1621,1542)(1621,762)(1621,762)
\path(1771,1542)(1771,732)
\path(1021,522)(1021,12)(1171,12)(1171,522)
\path(1621,522)(1621,12)(1771,12)(1771,522)
\path(1006,1842)(1006,3342)(1156,3342)(1156,1842)
\path(406,1842)(406,3342)(556,3342)(556,1842)
\path(631,2067)(931,2067)
\path(1231,2067)(1531,2067)
\path(406,552)(406,42)(556,42)(556,552)
\path(406,1587)(406,837)
\path(556,1617)(556,867)
\spline(331,2082)
(91,2082)(1,1932)
	(31,1752)(121,1692)(2101,1692)
	(2191,1782)(2221,1902)(2191,2022)
	(2071,2082)(1981,2082)
\path(1531,867)	(1484.009,878.602)
	(1456.000,867.000)

\path(1456,867)	(1421.198,814.443)
	(1409.598,754.500)
	(1421.198,694.557)
	(1456.000,642.000)

\path(1456,642)	(1505.790,601.395)
	(1561.117,572.391)
	(1620.136,554.988)
	(1681.000,549.188)
	(1741.864,554.988)
	(1800.883,572.391)
	(1856.210,601.395)
	(1906.000,642.000)

\path(1906,642)	(1940.802,694.557)
	(1952.402,754.500)
	(1940.802,814.443)
	(1906.000,867.000)

\path(1906,867)	(1877.991,878.602)
	(1831.000,867.000)

\path(931,867)	(884.009,878.602)
	(856.000,867.000)

\path(856,867)	(821.198,814.443)
	(809.598,754.500)
	(821.198,694.557)
	(856.000,642.000)

\path(856,642)	(905.790,601.395)
	(961.117,572.391)
	(1020.136,554.988)
	(1081.000,549.188)
	(1141.864,554.988)
	(1200.883,572.391)
	(1256.210,601.395)
	(1306.000,642.000)

\path(1306,642)	(1340.802,694.557)
	(1352.402,754.500)
	(1340.802,814.443)
	(1306.000,867.000)

\path(1306,867)	(1277.991,878.602)
	(1231.000,867.000)

\path(331,867)	(284.009,878.602)
	(256.000,867.000)

\path(256,867)	(221.198,814.443)
	(209.597,754.500)
	(221.198,694.557)
	(256.000,642.000)

\path(256,642)	(305.790,601.395)
	(361.117,572.391)
	(420.136,554.988)
	(481.000,549.188)
	(541.864,554.988)
	(600.883,572.391)
	(656.210,601.395)
	(706.000,642.000)

\path(706,642)	(740.802,694.557)
	(752.402,754.500)
	(740.802,814.443)
	(706.000,867.000)

\path(706,867)	(677.991,878.602)
	(631.000,867.000)

\put(2581,1767){\makebox(0,0)[lb]{$=$}}
\end{picture} }
        	\hspace{-1.9mm}
	\end{array}
 $$
\caption{The lantern identity in $P_3$.
Horizontal circles in the picture correspond (by a $+1$, i.e., a right-handed
twist)
framed pure braids on $3$-strands denoted by $\tau_i, \tau_{ij}, 
\tau_{ijk}$.}\lbl{lantern3}
\end{figure}

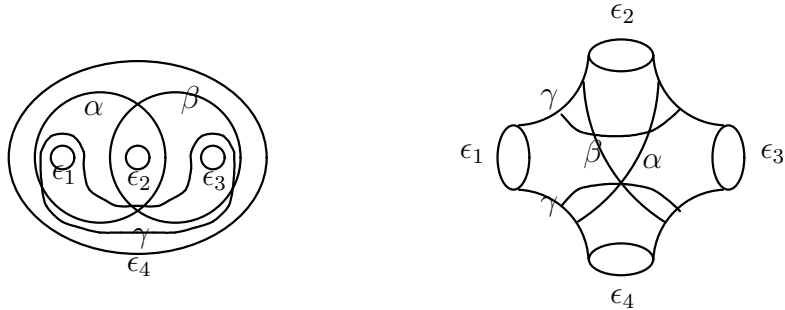
\begin{figure}[htpb]
$$ \printname{lantern2}
	\setlength{\unitlength}{0.03\standardunitlength}
	\begin{array}{c}  \hspace{-1.7mm}
        	\raisebox{-8pt}{\begingroup\makeatletter\ifx\SetFigFont\undefined
\def\x#1#2#3#4#5#6#7\relax{\def\x{#1#2#3#4#5#6}}%
\expandafter\x\fmtname xxxxxx\relax \def\y{splain}%
\ifx\x\y   
\gdef\SetFigFont#1#2#3{%
  \ifnum #1<17\tiny\else \ifnum #1<20\small\else
  \ifnum #1<24\normalsize\else \ifnum #1<29\large\else
  \ifnum #1<34\Large\else \ifnum #1<41\LARGE\else
     \huge\fi\fi\fi\fi\fi\fi
  \csname #3\endcsname}%
\else
\gdef\SetFigFont#1#2#3{\begingroup
  \count@#1\relax \ifnum 25<\count@\count@25\fi
  \def\x{\endgroup\@setsize\SetFigFont{#2pt}}%
  \expandafter\x
    \csname \romannumeral\the\count@ pt\expandafter\endcsname
    \csname @\romannumeral\the\count@ pt\endcsname
  \csname #3\endcsname}%
\fi
\fi\endgroup
\begin{picture}(10877,4155)(0,-10)
\thicklines
\put(7080.500,3610.500){\arc{2057.681}{0.0511}{1.5197}}
\put(10035.500,3610.500){\arc{2057.681}{1.6218}{3.0905}}
\put(7081.000,656.000){\arc{2056.631}{4.7630}{6.2326}}
\put(10035.000,656.000){\arc{2056.631}{3.1922}{4.6618}}
\put(6659.136,3281.864){\arc{4851.759}{0.0408}{1.0058}}
\put(10456.863,3281.864){\arc{4851.757}{2.1358}{3.1008}}
\put(8558,708){\ellipse{900}{450}}
\put(8558,3558){\ellipse{900}{450}}
\put(7058,2133){\ellipse{450}{900}}
\put(10058,2133){\ellipse{450}{900}}
\put(1808,2133){\ellipse{3600}{2700}}
\put(758,2133){\ellipse{336}{336}}
\put(1808,2133){\ellipse{336}{336}}
\put(2858,2133){\ellipse{336}{336}}
\put(1283,2133){\ellipse{1824}{1824}}
\put(2333,2133){\ellipse{1824}{1824}}
\path(7733,2733)	(7773.760,2684.172)
	(7809.591,2642.577)
	(7869.106,2578.451)
	(7916.818,2535.350)
	(7958.000,2508.000)

\path(7958,2508)	(7994.854,2492.107)
	(8040.022,2478.361)
	(8090.778,2466.620)
	(8144.398,2456.741)
	(8198.155,2448.582)
	(8249.324,2442.000)
	(8295.181,2436.854)
	(8333.000,2433.000)

\path(8333,2433)	(8399.228,2429.169)
	(8440.017,2428.212)
	(8483.000,2427.892)
	(8525.983,2428.212)
	(8566.772,2429.169)
	(8633.000,2433.000)

\path(8633,2433)	(8670.838,2437.172)
	(8716.743,2443.085)
	(8767.972,2450.612)
	(8821.782,2459.625)
	(8875.433,2469.998)
	(8926.181,2481.603)
	(8971.284,2494.313)
	(9008.000,2508.000)

\path(9008,2508)	(9080.612,2550.830)
	(9122.122,2580.586)
	(9164.645,2613.244)
	(9206.337,2646.808)
	(9245.355,2679.285)
	(9308.000,2733.000)

\path(9308,2733)	(9334.734,2758.369)
	(9383.000,2808.000)

\path(7733,1458)	(7773.760,1506.828)
	(7809.591,1548.423)
	(7869.106,1612.549)
	(7916.818,1655.650)
	(7958.000,1683.000)

\path(7958,1683)	(7994.854,1698.893)
	(8040.022,1712.639)
	(8090.778,1724.380)
	(8144.398,1734.259)
	(8198.155,1742.418)
	(8249.324,1749.000)
	(8295.181,1754.146)
	(8333.000,1758.000)

\path(8333,1758)	(8399.228,1761.831)
	(8440.017,1762.788)
	(8483.000,1763.108)
	(8525.983,1762.788)
	(8566.772,1761.831)
	(8633.000,1758.000)

\path(8633,1758)	(8670.838,1753.828)
	(8716.743,1747.915)
	(8767.972,1740.388)
	(8821.782,1731.375)
	(8875.433,1721.002)
	(8926.181,1709.397)
	(8971.284,1696.687)
	(9008.000,1683.000)

\path(9008,1683)	(9080.612,1640.170)
	(9122.122,1610.414)
	(9164.645,1577.756)
	(9206.337,1544.192)
	(9245.355,1511.715)
	(9308.000,1458.000)

\path(9308,1458)	(9334.734,1432.631)
	(9383.000,1383.000)

\path(458,2208)	(479.561,2256.920)
	(499.440,2298.581)
	(535.910,2362.759)
	(608.000,2433.000)

\path(608,2433)	(676.663,2456.321)
	(716.540,2462.152)
	(758.000,2464.095)
	(799.460,2462.152)
	(839.337,2456.321)
	(908.000,2433.000)

\path(908,2433)	(958.127,2391.229)
	(1001.799,2330.888)
	(1036.071,2265.352)
	(1058.000,2208.000)

\path(1058,2208)	(1064.017,2157.286)
	(1061.754,2094.619)
	(1057.613,2032.392)
	(1058.000,1983.000)

\path(1058,1983)	(1066.781,1915.509)
	(1073.143,1874.448)
	(1081.096,1831.410)
	(1090.850,1788.571)
	(1102.613,1748.106)
	(1133.000,1683.000)

\path(1133,1683)	(1178.031,1639.238)
	(1240.895,1597.699)
	(1306.062,1561.310)
	(1358.000,1533.000)

\path(1358,1533)	(1430.608,1491.776)
	(1472.533,1470.949)
	(1508.000,1458.000)

\path(1508,1458)	(1553.109,1450.946)
	(1613.457,1448.595)
	(1651.544,1449.183)
	(1696.077,1450.946)
	(1747.936,1453.885)
	(1808.000,1458.000)

\path(458,2208)	(456.602,2147.408)
	(455.604,2095.172)
	(455.005,2050.412)
	(454.805,2012.250)
	(455.604,1952.203)
	(458.000,1908.000)

\path(458,1908)	(461.528,1862.908)
	(465.906,1808.123)
	(471.545,1746.896)
	(478.854,1682.483)
	(488.243,1618.133)
	(500.124,1557.101)
	(514.906,1502.639)
	(533.000,1458.000)

\path(533,1458)	(574.529,1398.225)
	(634.752,1333.759)
	(700.349,1275.164)
	(758.000,1233.000)

\path(758,1233)	(806.490,1210.651)
	(869.360,1189.834)
	(932.800,1171.849)
	(983.000,1158.000)

\path(983,1158)	(1048.452,1139.190)
	(1088.931,1127.909)
	(1131.718,1116.386)
	(1174.664,1105.378)
	(1215.625,1095.640)
	(1283.000,1083.000)

\path(1283,1083)	(1349.399,1079.414)
	(1390.251,1079.265)
	(1433.270,1079.812)
	(1476.255,1080.759)
	(1517.006,1081.805)
	(1583.000,1083.000)

\path(1583,1083)	(1616.048,1083.000)
	(1661.049,1083.000)
	(1723.275,1083.000)
	(1762.495,1083.000)
	(1808.000,1083.000)

\path(3158,2208)	(3136.439,2256.920)
	(3116.560,2298.581)
	(3080.090,2362.759)
	(3008.000,2433.000)

\path(3008,2433)	(2939.337,2456.321)
	(2899.460,2462.152)
	(2858.000,2464.095)
	(2816.540,2462.152)
	(2776.663,2456.321)
	(2708.000,2433.000)

\path(2708,2433)	(2657.873,2391.229)
	(2614.201,2330.888)
	(2579.929,2265.352)
	(2558.000,2208.000)

\path(2558,2208)	(2551.983,2157.286)
	(2554.246,2094.619)
	(2558.387,2032.392)
	(2558.000,1983.000)

\path(2558,1983)	(2549.219,1915.509)
	(2542.857,1874.448)
	(2534.904,1831.410)
	(2525.150,1788.571)
	(2513.387,1748.106)
	(2483.000,1683.000)

\path(2483,1683)	(2437.969,1639.238)
	(2375.105,1597.699)
	(2309.938,1561.310)
	(2258.000,1533.000)

\path(2258,1533)	(2185.392,1491.776)
	(2143.467,1470.949)
	(2108.000,1458.000)

\path(2108,1458)	(2062.891,1450.946)
	(2002.542,1448.595)
	(1964.456,1449.183)
	(1919.923,1450.946)
	(1868.064,1453.885)
	(1808.000,1458.000)

\path(3158,2208)	(3159.398,2147.408)
	(3160.396,2095.172)
	(3160.995,2050.412)
	(3161.195,2012.250)
	(3160.396,1952.203)
	(3158.000,1908.000)

\path(3158,1908)	(3154.472,1862.908)
	(3150.094,1808.123)
	(3144.455,1746.896)
	(3137.146,1682.483)
	(3127.757,1618.133)
	(3115.876,1557.101)
	(3101.094,1502.639)
	(3083.000,1458.000)

\path(3083,1458)	(3041.471,1398.225)
	(2981.247,1333.759)
	(2915.651,1275.164)
	(2858.000,1233.000)

\path(2858,1233)	(2809.510,1210.651)
	(2746.640,1189.834)
	(2683.200,1171.849)
	(2633.000,1158.000)

\path(2633,1158)	(2567.548,1139.190)
	(2527.069,1127.909)
	(2484.283,1116.386)
	(2441.336,1105.378)
	(2400.375,1095.640)
	(2333.000,1083.000)

\path(2333,1083)	(2266.601,1079.414)
	(2225.749,1079.265)
	(2182.730,1079.812)
	(2139.745,1080.759)
	(2098.994,1081.805)
	(2033.000,1083.000)

\path(2033,1083)	(1999.952,1083.000)
	(1954.951,1083.000)
	(1892.725,1083.000)
	(1853.505,1083.000)
	(1808.000,1083.000)

\put(8408,4008){\makebox(0,0)[lb]{$\e_2$}}
\put(6308,2058){\makebox(0,0)[lb]{$\e_1$}}
\put(10508,2058){\makebox(0,0)[lb]{$\e_3$}}
\put(8408,33){\makebox(0,0)[lb]{$\e_4$}}
\put(8033,1983){\makebox(0,0)[lb]{$\b$}}
\put(8858,1983){\makebox(0,0)[lb]{$\a$}}
\put(7433,2808){\makebox(0,0)[lb]{$\g$}}
\put(7433,1308){\makebox(0,0)[lb]{$\g$}}
\put(2408,2733){\makebox(0,0)[lb]{$\b$}}
\put(1058,2733){\makebox(0,0)[lb]{$\a$}}
\put(608,1758){\makebox(0,0)[lb]{$\e_1$}}
\put(1658,1683){\makebox(0,0)[lb]{$\e_2$}}
\put(2708,1683){\makebox(0,0)[lb]{$\e_3$}}
\put(1733,858){\makebox(0,0)[lb]{$\g$}}
\put(1658,483){\makebox(0,0)[lb]{$\e_4$}}
\end{picture} }
        	\hspace{-1.9mm}
	\end{array}
 $$
\caption{Two views of the curves that appear in the 
lantern identity on a sphere.}\lbl{lantern2}
\end{figure}
   
\section{The \LT\ and \fti s}
\subsection{Finite type invariants and the \LT}

Using the following elementary lemma:
\begin{lemma}
\lbl{lem.elem}
If $R=\BZ \la a,b \ra$ is the noncommutative ring on $\{ a, b \}$
 and $I \eqdef (\ov a, \ov b)$ the two-sided augmentation ideal 
(where $\ov x \eqdef 1-x$), then we have:
\begin{equation}
\ov {ab} \equiv \ov a + \ov b \text{ mod } I^2
\end{equation}
\end{lemma}
\noindent
together with the lantern identity in $P_3$, we deduce that:
\begin{proposition}\cite{Oh}
\lbl{cor.3bands}
For $R=\BZ P_3$ and $I$ the augmentation ideal, we have:
\begin{equation}
\lbl{eq.3bands}
\ov \t_{123} \equiv  \ov \t_{12} + \ov \t_{13} + \ov \t_{23} - \ov \t_{1} 
- \ov \t_{2} -  \ov \t_{3} \text{ mod } I^2
\end{equation} 
\end{proposition}
Together with the moves of  \cite[Figure 2.1]{Oh},
Ohtsuki used the above equation to show that the vector space
of finite type $m$ invariants of integral homology 3-spheres is finite
dimensional for every nonnegative integer $m$, \cite{Oh}. 
It was shown in \cite{L} (see also \cite[Remark 2.9]{GL})
that the moves of \cite[Figure 2.1]{Oh} generate the surgery equivalence
relation of algebraically split links (i.e., links with linking number zero); 
however  a conceptual understanding of 
 equation \eqref{eq.3bands} was missing. The \LT\ provides such an
explanation.

A few questions are in order:
\begin{question}
Can we improve Lemma \ref{lem.elem} and Corollary \ref{cor.3bands}
mod $I^m$ for any $ m > 2$? 
\end{question}
 The answer is positive, as follows.
\begin{lemma}
\lbl{lem.elem2}
If $R=\BZ \la a, b, b^{-1} \ra $, then  in the completion $\hat{R}_I$
with respect to the augmentation ideal, we have: 
$$
1- \ov {ab^{-1}} = (1- \ov a)(1- \ov b)^{-1}
$$
\end{lemma}
\begin{proof}
Since $a= 1- \ov a  $ and $b^{-1}=1/b=(1-\ov b)^{-1}$, it follows
that $1-\ov {ab^{-1}}= ab^{-1}=(1- \ov a)(1- \ov b)^{-1}$.
\end{proof}

On the other hand, the lantern identity implies that for any number 
$n \geq 3$ of strands
we have in $P_n$:
$$
\t_{12 \dots n} \prod_{i=1}^n \t_i^{n-2} = \prod_{1 \leq i < j \leq n}
\t_{ij}
$$
(where $\t_i$ commute with $\t_{ij}$, and the procust is taken
lexicographically). Thus, we deduce 
that in the group ring $\widehat{\BZ P_n}$ (completed with respect to
the augemntation ideal) we have:
$$
1-\ov \t_{12 \dots n}= \frac{\prod_{1 \leq i < j \leq n}
(1 - \ov \t_{ij})}{\prod_{i=1}^n(1 - \ov \t_i)^{n-2}}
$$
(where the product is taken lexicographically). This was first
obtained in \cite[Theorem 4]{GL}. The reader is invited to
be convinced that the above identity is indeed a generalization
of Proposition \ref{cor.3bands}.

\begin{question}
Do we really need the \LT\ in order to show that
the space of finite type $m$ invariants of \ihs s
is finite dimensional? 
\end{question}

Ohtsuki's proof uses only the fact that we can express $\ov \tau_{123}$
of equation \eqref{eq.3bands} in terms of $\ov \tau_{i}, \ov \tau_{ij}$. 
Therefore any fictional identity on $P_3$ that expressed a nonzero
power of $\tau_{123}$ in terms of $\tau_{i}, \tau_{ij}$ would suffice.

\subsection{Finite $T$-type invariants and the \LT }
 
We begin by recalling a few definitions and notation from \cite{GL2}.
All manifolds will be oriented and all maps will be orientation preserving.
Let $\M$ denote the vector space (over $\BQ$) on the set of (oriented)
\ihs s.
Given an  embedding $\S_g \to M$ of a closed genus
$g$ surface in an \ihs\ $M$, let $\T_g$ denote the Torelli group of
$\S_g$ (i.e., the group of (orientation preserving) diffeomorphism classes
of surface diffeomorphisms that act trivially in the homology), let
$\BQ \T_g$ denote the group-ring and $I$ its augmentation ideal. 
The process of cutting $M$ across $\S_g$, twisting by an element of $\T_g$
and gluing back, defines (by linear extension) a map:
$\Phi_f: \BQ \T_g \to \M$. In \cite[Definition 1.1]{GL2} we considered the
decreasing filtration $\FT \ast$ on $\M$ defined by $\FT m= 
\cup_{\text{all } f} \Phi_f(I^m)$. We call a map $\lambda: \M \to \BQ$
an invariant of \ihs s of $T$-{\em type} $m$ if $\lambda(\FT {m+1})=0$.
As an application of the \LT\ we show that:
\begin{proposition}
The vector space of $T$-type $m$ invariants of \ihs s is finite dimensional
for every nonnegative integer $m$. Furthermore, the graded vector space
 is zero dimensional for odd $m$.
\end{proposition}

\begin{proof}
It suffices to show that  the graded vector space
$\GT m \eqdef \FT m/\FT {m+1}$ is finite (resp. zero) dimensional
for every nonnegative (resp. odd) integer $m$. A geometric argument
of \cite[Proposition 1.16]{GL2} implies that $\FT m$ is the union over all 
Heegaard surfaces $h: \S_g \to M$ in \ihs s $M$, and thus $\GT m=\cup_{h} 
\Phi_h(I^m/I^{m+1})$. The action of the mapping
class group $\Ga_g$ on $\T_g$ by conjugation implies that  
that $I^m/I^{m+1}$ is a $Sp(H)$ (and thus, a
$GL(L^+)$) module, where $H=H_1(\S_g, \BZ)$,
$L^{\pm}= \text{Ker}(h: H \to H_1(M^{\pm},\BZ))$ and $M= M^+ \cup_h
M^-$.  Furthermore, another geometric argument of 
\cite[Lemmas 3.1-3.3]{GL3} implies that $\cup_h \Phi_h(I^m/I^{m+1})=  
\cup_{h'}\Phi_h((I^m/I^{m+1})^{GL(L^+)})$, for all $h'$ standard Heegaard
splittings of $S^3$ (of arbitrary genus). 

D. Johnson \cite{Jo2} introduced a group homomorphism (well known
as the Johnson homomorphism) $\T_g \to U\eqdef\La^3H/H$. 
Using the lantern identity, he 
showed in \cite{Jo5} that the Johnson homomorphism coincides (rationally) 
with the
 abelianization of the Torelli group $\T_g$ as a $Sp(H)$ module.

Furthermore, both the Johnson homomorphism and the map $\Phi_{h'}$ above
are  {\em stable} with respect to the inclusion of a surface to 
another, in the following sense: for a surface $\S$ with one boundary 
component included in a 
 closed surface $\S'$, let $\hat{\S}$ denote the
surface obtained by closing $\S$ by the addition of a disc. Then, we have
a canonical inclusion $H_1(\hat{\S},\BZ) \simeq H_1(\S, \BZ) \to
H_1(\S', \BZ)$, and both the Johnson homomorphism and the map $\Phi_{h'}$
respect this inclusion. 

Since $I/I^2$ can also be (rationally) identified with the abelianization
of $\T_g$, we deduce that $I/I^2 \simeq U$ as a $Sp(H)$ module. Thus,
for a fixed $m$, $I^m/I^{m+1}$ is a quotient of $\otimes^mU$, and hence
of $\La^3 H$. Since 
$\La^3 H= \La^3 L + \La^2 L \otimes L^\star + \La^2 L^\star \otimes L +
 \La^3 L^\star$ (where $L=L^+$ and $L^\ast$ is the dual of $L$) is a
$GL(L^+)$ representation with odd heighest weights odd of length at most $3$, 
if follows from classical invariant theory \cite{W} 
and the above mentioned stability of the Johnson homomorphism
with respect to $g$ that  
the vector spaces $(I^m/I^{m+1})^{GL(L^+)}$  are naturally isomorphic,
 finite dimensional for genus $g \geq 3m$ and in addition zero dimensional
for odd $m$.  
Furthermore, by the stability of $\Phi_h$ we have that
 $\GT m= \Phi_{h_{3m}}((I^m/I^{m+1})^{GL(L^+)})$
where $h_{3m}$ is the standard genus $3m$ Heegaard splitting of $S^3$.  
This finishes the proof of the proposition.
\end{proof}

\begin{remark}
In \cite[Theorem 4]{GL2} we have shown among other things
that $T$-type $2m$ invariants of \ihs s coincide with $T$-type $2m-1$
and with type $3m$ invariants of \ihs s. 
Thus, the above proposition follows from the fact type $m$ invariants of
\ihs s form a finite dimensional vector space. Compared to the
rather involved proof of \cite[Theorem 4]{GL2}, the  above proof is
new and shorter; furthermore it  reduces the finiteness result
of the proposition to the fact that the abelianization of the Torelli group
is finitely generated.
\end{remark}

It seems natural to ask the following:
\begin{question}
Do we really need the \LT\ in order to show that
the space of $T$-type $m$ invariants of \ihs s
is finite dimensional? 
\end{question}

The above proof shows that as long as (rationally) the abelianization of the 
Torelli is finitely generated, $T$-type $m$ invariants of \ihs s form
a finite dimensional space. In addition, as long as the abelianization of the
Torelli group  as a  $Sp(H)$ module is  included in  an odd tensor
power of $H$, the graded vector space of finite $T$-type $m$ invariants
is zero dimensional for odd $m$. The ambitious reader 
may figure out a  form of a fictional \lt\   compatible with
Johnson's arguments which show that the Torelli group (and thus
its abelianization) is finitely
generated \cite{Jo4} or with Johnson's arguments that 
 determine the abelianization of the Torelli group \cite{Jo5}.

\section{Just for completeness}
\lbl{sec.comp}

Just for completeness, in this section we mention two more applications
of the \lt\ to the following finiteness theorems:
\begin{itemize}
\item
The abelianization of the mapping class group of a closed genus $g > 2$
surface is trivial, due to J. Powell \cite{Po} and independently due to
J. Harer \cite{Ha}. 
\item
The conformal dimensions and the central charge of an arbitrary topological
quantum field theory in 3 dimensions are rational numbers,
due to C. Vafa \cite{Va}.
\end{itemize}

We summarize Harer's argument here. Since the mapping class group is generated
by Dehn twists on simple closed curves, the \lt\ implies that it is generated
by twists on nonseparating curves. Since any two nonseparating
simple closed curves are conjugate, the lantern identity
(applied so that the  Dehn twists on all
seven participating curves are nonseparating and right-handed) 
implies the abelianization of the the mapping class group is trivial.

For a reproduction of Vafa's argument as well as a friendly and complete
 definition of the terms and concepts involved, see \cite{BK}.

The ambitious reader may figure out a form of a fictional \lt\ which 
would still make J. Harer's and C. Vafa's arguments  work.

\subsection{Acknowledgment}

We wish to thank Jerome  Levine for encouraging conversations.


\ifx\undefined\bysame
	\newcommand{\bysame}{\leavevmode\hbox to3em{\hrulefill}\,}
\fi


\begin{thebibliography}{[EMSS]}

\bibitem[BK]{BK} B. Bakalov, A. Kirillov Jr.,
       {\em Tensor categories}, to appear.

\bibitem[GL1]{GL} S. Garoufalidis, J. Levine,
       {\em On finite type 3-manifold invariants II},  
       Math. Annalen, {\bf 306} (1996) 691--718.

\bibitem[GL2]{GL2} \bysame,
       {\em Finite type 3-manifold invariants, the mapping 
       class group and blinks}, Journal of Diff. Geom., in press.

\bibitem[GL3]{GL3} \bysame,
        {\em Finite type 3-manifold invariants and the Torelli group I},
        Inventiones, in press.

\bibitem[Ha]{Ha} J. Harer,
       {\em The second homology of the mapping class group of an orientable
       surface}, Inventiones {\bf 72} (1983) 221--239.

\bibitem[Jo1]{Jo1} D. Johnson, 
       {\em Homeomorphisms of a surface which act trivially on homology},
       Proc. Amer. Math. Soc., {\bf 75} (1979) 118--125.

\bibitem[Jo2]{Jo2}\bysame,
       {\em An abelian quotient of the mapping class group $\mathcal T_g$},
       Math. Annalen {\bf 249} (1980) 225--242.

\bibitem[Jo3]{Jo3}\bysame ,
       {\em A survey of the Torelli group},
       Contemporary Math. {\bf 20} (1983) 163--179.

\bibitem[Jo4]{Jo4}\bysame,
       {\em The structure of the Torelli group I: a finite set of generators
       for $\mathcal T_g$},
       Annals Math. {\bf 118} (1983) 423--442.

\bibitem[Jo5]{Jo5}\bysame,
       {\em The structure of the Torelli group III: the abelianization
       of $\mathcal T_g$},
       Topology, {\bf 24} (1985) 127--144.

\bibitem[L]{L} J. Levine, {\em Surgery equivalence of links}, 
       Topology {\bf 26} (1987) 45--61.

\bibitem[Oh]{Oh} T. Ohtsuki,
        {\em Finite type invariants of integral homology 3-spheres}, 
        J. Knot Theory and its Rami. {\bf 5} (1996) 101--115. 

\bibitem[Po]{Po} J. Powell,
        {\em Two theorems on the mapping class group of surfaces},
        Proc. Amer. Math. Soc., {\bf 68} (1978) 347--350.

\bibitem[Va]{Va} C. Vafa,
        {\em Towards classification of conformal field theories},
        Phys. Lett. B {\bf 206} (1988) 421--426.

\bibitem[W]{W} H.Weyl,
	{\em The classical groups, their invariants and representations},
	2nd ed., Princeton Univ. Press, Princeton NJ, 1946.

\end{thebibliography}
\end{document}